\documentclass[12pt]{article}

\usepackage{amsgen,amsmath,amstext,amsbsy,amsopn,amsfonts,amssymb}
\newtheorem{theorem}{Theorem}
\newtheorem{lemma}[theorem]{Lemma}
\newtheorem{corollary}[theorem]{Corollary}
\newtheorem{proposition}[theorem]{Proposition}
 \newtheorem{defi}[theorem]{Definition}

\newtheorem{exa}[theorem]{Example}

\newtheorem{rem}[theorem]{Remark}
\newenvironment{remark}{\begin{rem}\rm}{\end{rem}}
\newtheorem{rems}[theorem]{Remarks}

\newtheorem{ack}[theorem]{Acknowlegment}

\def\bsq{\blacksquare\medskip}
\def\n{\noindent}

\let \g =\gamma
\let \n = \noindent

\let \f =\varphi

\def\Ker{\mathop{\rm Ker}}

\def\fV {\varphi\big\vert_{V}}

\def\H{\mathcal H}
\def\L{\mathcal L}
\def\K{\mathcal K}

\def\O{\mathcal O}
\def\CalS{\mathcal S}
\def\Hpi{(\pi, {\mathcal H})}

\def\NN{{\mathbb N}}
\def\ZZ{{\mathbb Z}}
\def\CCC{{\mathbb C}}
\def\RRR{{\mathbb R}}
\def\QQ{\mathbb Q}
\def\RR+{{\mathbb R}^*}
\def\TT{\mathbb T}
\def\KK{\mathbb K}
\def\FF{\mathbb F}

\def\Q_p{{\mathbb Q}_p}
\def\S1{{\mathbb S}^1}

\def\eps{\varepsilon}
\def\Ga{\Gamma}
\def\ga{\gamma}
\def\vfi{\varphi}

\def\SL{SL_n(\ZZ)}
\def\Sl{SL_{n-1}(\ZZ)}

\def\PSL{PSL_n(\ZZ)}
\def\tous{\qquad\text{for all}\quad}
\def\et{\qquad\text{and}\qquad}
\begin{document}

\title{Operator-algebraic superridigity  for $\SL,\ n\geq 3$}
\author{Bachir Bekka}
\date{\today }

\maketitle

\begin{abstract}
For $n\geq 3,$ let $\Ga=\SL.$ 
We prove the following superridigity result for $\Ga$
in the context of operator algebras. 
Let $L(\Ga)$ be the von Neumann algebra
generated by the left regular representation
of $\Ga.$
Let $M$ be a finite factor and let 
$U(M)$ be its unitary group. Let $\pi: \Ga\to U(M)$ be a group
homomorphism such that $\pi(\Ga)''=M.$
Then either
\begin{itemize}
\item[(i)] $M$ is finite dimensional, or 
\item [(ii)] there exists a subgroup
of finite index $\Lambda$ of $\Ga$ such 
that $\pi|_\Lambda$ extends to a homomorphism 
$U(L(\Lambda))\to U(M).$ 
\end{itemize}
This answers, in the special case of $\SL,$ 
a question of A.~Connes discussed in \cite[Page 86]{Jones}.
The result is deduced from  a complete description 
of the tracial states on  the full
$C^*$--algebra of $\Ga.$

As another application,  we show that
the full $C^*$--algebra of
 $ \Ga$ has no faithful tracial state, thus
answering a question of E.~Kirchberg.

\end{abstract}

\section{Introduction}
Two  major achievements in the study of
discrete subgroups in semi-simple Lie groups
  are Mostow's ridigity theorem
and Margulis' superrigidity theorem.
A weak version of the latter is as
follows. Let $\Ga$ be a lattice
in a simple real Lie group $G$
with finite centre and
with $\RRR-{\rm rank} (G) \geq 2.$
Let $H$ be another  simple 
real Lie group with finite centre,
and let
$\pi:\Ga\to H$ be a homomorphism
such that $\pi(\Ga)$ is Zariski-dense in $H.$
  Then, either $H$ is compact or 
there exists a finite index subgroup
$\Lambda$ of $\Ga$ such that $\pi\vert_\Lambda$ extends to a 
continuous homomorphism $G\to H.$
For more general results, see \cite{Margulis}  and \cite{Zimmer}. Moreover, 
as shown by Corlette,  the superridity theorem
continues to hold for the  simple real Lie groups $G$
with $\RRR-{\rm rank} (G) =1$ which are not locally isomorphic
to $SO(n,1)$ or $SU(n,1)$.

In the theory of von Neumann algebras, 
discrete groups (as well as their actions) 
always played a prominent r\^ole.
To  a discrete group $\Ga$ 
is associated a distinguished von  Neumann algebra
$L(\Ga),$ namely
the von Neumann algebra generated by the left regular representation
$\lambda_\Ga$ of $\Ga;$ thus,
$L(\Ga)$ is the closure for the strong operator topology
of the
linear span of $\{\lambda_\Ga(\ga)\ :\ \ga\in\Ga\}$
in the algebra $\L(\ell^2(\Ga))$ 
of all bounded operators on the Hilbert space
$\ell^2(\Ga).$

The first rigidity result in the context of operator algebras
is  the result by A. Connes \cite{Conne--80}
showing that, for   a  group $\Ga$ with Kazhdan's Property (T),
  the group of outer automorphisms of  $L(\Gamma)$  
is countable. A major  problem in this area
is whether  such a group $\Ga$
can be reconstructed from its  von Neumann
algebra $L(\Gamma).$ 
 In recent years, a series of 
remarkable results concerning this question,
with applications to ergodic theory, have been obtained
by S.~Popa (\cite{Popa1}, \cite{Popa2}; 
for an account, see \cite{Vaes}).
Other relevant work 
includes \cite{Cowling-Haagerup} and \cite{Furman}.

The purpose of this paper is to discuss
another kind of rigidity, namely the rigidity
of  a discrete group in the unitary group
of its von Neumann algebra.
If $\Ga$ is a discrete group,
we  view $\Ga$ as a subgroup of 
the unitary group $U(L(\Ga))$ of $L(\Ga),$ that
is, the group of  the unitary operators in $L(\Ga).$ 
 It was suggested by Connes (see \cite[Page 86]{Jones})
that, for $\Ga$ as in the statement of Margulis'
theorem,  a superrigity result should hold
in which $G$ above is replaced
by $U(L(\Ga))$ and $H$ by 
the unitary group $U(M)$ of  a type $II_1$
factor.
We prove such a superrigity result
in the case $\Ga= SL_n(\ZZ)$ for $n\geq 3.$

Recall that a von Neumann algebra $M$  is a factor
if the centre of $M$ is reduced to the scalar operators.
The von Neumann algebra $M$ is said to be finite
if there exists a finite
normal faithful trace  on $M.$ 
A finite factor  is a type $II_1$ factor which
is infinite dimensional.
Recall also that $L(\Ga)$ is a finite von Neumann
algebra. Moreover, $L(\Ga)$ is a factor if and only if
$\Ga$ is an ICC-group, that is, if all  its conjugacy classes, except $\{e\},$ are infinite.
For an account on the theory of von Neumann algebras,
see \cite{Dixmier-vN}.

\begin{theorem}
\label{Theo2}
 Let $\Ga=SL_n(\ZZ)$ for $n\geq 3.$
Let $M$ be a  finite factor and
let $U(M)$ its unitary group. Let 
$\pi:\Ga\to U(M)$ be a group homomorphism.
Assume that the linear span of $\pi(\Ga)$ 
in dense in $M$ for the strong operator topology.
Then either
\begin{itemize}
\item[(i)] $M$ is finite dimensional, that is, $M$ is isomorphic to
a matrix algebra $M_n(\CCC)$ for some $n\in\NN,$ in which case
 $\pi$ factorizes to  a multiple of an irreducible 
 representation of some congruence quotient
$SL_n({\ZZ}/N{\ZZ})$ for $N\in\NN$, or 
\item [(ii)] there exists a subgroup
of finite index $\Lambda$ of $\Ga$ such 
that $\pi|_\Lambda$ extends to a normal homomorphism 
$L(\Lambda)\to M$  of von Neumann algebras.
In particular, $\pi|_\Lambda$ extends to a group
homomorphism $U(L(\Lambda))\to U(M).$
\end{itemize}
\end{theorem}
Let $C$ be the centre of $SL_n(\ZZ).$
Observe that $C$ is trivial for odd $n$ 
and $C=\{\pm I\}$ for even  $n.$ 
If, in the statement of the theorem above, we 
take instead $\Ga= PSL_n(\ZZ)= SL_n(\ZZ)/C,$
then $L(\Ga)$ is a factor 
and the conclusion (ii) holds for $\Lambda= \Gamma.$

The method of proof of Theorem~\ref{Theo2} can be adapted 
to establish the same result when $\Ga$ is  the symplectic
group $Sp_{2n}(\ZZ)$ for $n\geq 2;$
it works presumably 
for the group of integral points of any Chevalley group of rank $\geq 2.$
No such result can be true for the modular group
$SL_2(\ZZ);$ see Remark~\ref{Rem-Free group} below.

\begin{remark} 
\label{Rem-Guyan} (i) Let $\Ga$ be a 
countable ICC--group with Kazhdan's Property
(T). It was shown in \cite{CoJo} that 
$L(\Ga)$ cannot be a subfactor of $L(F_2)$,
where $F_2$ is a  non-abelian free group.
This, combined with Theorem~\ref{Theo2},
shows that every representation  of  $PSL_n(\ZZ)$ for $n\geq 3$
into $U(L(F_2))$
decomposes as  a direct sum of finite dimensional
representations.  This is a special case of
 a result of G.~Robertson \cite{Guyan}
valid for all groups with Property
(T).  For a related work, see \cite{Alain}.

\n
(ii) Let $M$ be a finite factor. 
 The unitary group $U(M)$ has as centre a copy
of the circle group $S^1,$  namely
the unitary scalar operators.
It was shown in \cite{Pierre} that
the projective unitary group $U(M)/S^1$
of $M$  is a simple group.
\end{remark}

The result in Theorem~\ref{Theo2} amounts to the
classification 
of the  \emph{characters} of $\Ga$ (see Section~\ref{S:Superigidite}),
that is, the functions $\f:\Ga\to \CCC$ with the following properties:
\begin{itemize}
\item $\f$ is central, that is, $\f(\ga x\ga^{-1})= \f(x)$ for all $\ga,x\in \Ga,$
\item $\f$ is positive definite, that is, $\sum_{i=1}^n \overline{c_j}c_i\f(\ga_j^{-1} \ga_i)\geq 0$
for all $\ga_1, \dots, \ga_n\in \Ga$ and $c_1,\dots, c_n\in \CCC,$ 
\item $\f$ is normalized, that is, $\f (e)=1,$
\item $\f$ is indecomposable, that is,  $\f$ cannot be written in a non-trivial
way as a convex combination of two central  positive definite normalized functions. 
\end{itemize}


There are two obvious examples
of characters of a group $\Ga.$
 First of all, the  normalized character (in the usual
sense) of an irreducible
 finite dimensional  unitary representations of $\Ga$
is a character of $\Ga$ in the above sense.
For $\Ga =\SL$, $n\geq 3$, it is well--known that every
such representation factorizes through some
congruence quotient $SL_n(\ZZ/N\ZZ)$ for an integer
$N;$ this a consequence of the solution of the congruence
subgroup problem (see \cite{BLS}, \cite{Mennicke}; see also \cite{Steinberg}). 
Much is known
about the characters of the finite groups $SL_n(\ZZ/N\ZZ)$; see \cite{Zel}.

Let $C$ be the centre
of the group $\Ga.$ Assume that all conjugacy classes, except those of the
elements from $C$, are infinite. Then,
for every unitary character $\chi$ of the abelian  group $C$, the trivial extension
$\widetilde\chi$ of $\chi$ to $\Ga,$ defined
by $\widetilde\chi=0$ on $\Ga\setminus C,$ 
is a character of $\Ga.$ 
In particular, if $\Ga$ is ICC, then $\delta_e,$ the Dirac function
at $e,$ is a character of $\Ga.$
When $n$ is even, all conjugacy classes of $PSL_n(\ZZ),$ except
$\{I\}$ and $\{-I\}$, are
infinite.

Our main result
says that  $\SL$ for 
$n\geq 3$  has no  characters
other than the obvious ones described above.

\par
\begin{theorem}
\label{Theo1}
Let $\varphi$ be a character of $\SL$ for 
$n\geq 3$. Then, either
\begin{itemize}
\item[(i)] $\varphi$ is the character
 of an irreducible  finite dimensional
 representation of some congruence quotient
$SL_n({\ZZ}/N{\ZZ})$ for $N\geq 1$, or 

\item[(ii)] $\varphi$ is  the trivial extension
$\widetilde\chi$ of a character $\chi$ of the centre of $\SL.$.
 \end{itemize}
\end{theorem}
\begin{remark} 
\label{Rem-Free group}
No classification of the characters of the modular group
$SL_2(\ZZ)$ can be expected. Indeed, this group 
contains the free non-abelian group $F_2$ 
on two generators as normal subgroup.
Every character of $F_2$ 
extends to a character on $SL_2(\ZZ).$ 
Now, $F_2$  
has a huge number of characters: if $M$ is any  finite factor
with  trace $\tau,$ 
every pair of unitaries in $M$ defines
a homomorphism $\pi:F_2\to U(M)$ and 
a corresponding character $\tau\circ\pi$ on $F_2.$
\end{remark}

\par
The problem of the description of the characters
of a discrete group $\Ga$ has been considered by
several authors.
 E.~Thoma \cite{ThomaSym} solved this problem for
 the infinite symmetric group $S_{\infty}$ (see also  \cite{Vershik}), 
 H~-~L.~Skudlarek \cite{Skud} for the group $\Ga=GL(\infty, \FF)$, where  $\FF$ 
is a finite field,
 and  D.~Voiculescu \cite{Voiculescu}  for $\Ga =U(\infty)$; see also \cite{StraVoi}
and \cite{Boyer}.
A.A.~Kirillov \cite{Kirillov} described the characters
of  $\Ga= GL_n( \KK)$ or  $ SL_n( \KK)$ for $n\geq 2,$ where $\KK$ is an infinite
field (see also \cite{Rosenberg} and  \cite{Ovc}).

Our proof of Theorem~\ref{Theo1} is based on an analysis
of the restriction $\f|_V$ of a given character $\f$ of $\SL$ to various
 copies $V$ of $\ZZ^{n-1}.$ 
We will see that we have a dichotomy
 corresponding to the two different types
of characters from Theorem~\ref{Theo1}:
either the measure on the torus $\TT^{n-1}$
associated to $\f|_V$ is atomic or this measure is the Lebesgue measure
for every  $V.$ 
An important ingredient in our analysis
is the solution of the congruence subgroup
for $\SL $ for $n\geq 3.$

The result of Theorem~\ref{Theo1} can be interpreted as a classification
of the traces on the full $C^*$-algebra $C^*(\Ga)$ of $\Ga=\SL$ for $n\geq 3$
(see Section~\ref{S:Factors}).

E.~Kirchberg asked in \cite[Remark 8.2, page 487]{Kirchberg} whether the full
$C^*$--algebra of $SL_4(\ZZ)$ has a faithful trace.
He was motivated by the fact that a positive answer
to this question would imply a series of  outstanding
conjectures in the theory of von Neumann algebras (see Section~\ref{S:Kirchberg}).
As a consequence of Theorem~\ref{Theo1}, we will see
that the answer to Kirchberg's question is negative, namely:
\par
\begin{corollary}
\label{Theo4}
The  full $C^*$--algebra of $\SL$ has no faithful tracial state
for $n\geq 3$.
\end{corollary}
In fact, we will prove
the stronger result Corollary~\ref{Cor-Kirch} below.

Recall that the reduced $C^*$-algebra
$C_r^*(\Ga)$ of a group $\Ga$ is the
closure of the linear span
of $\{\lambda_{\Ga}(\ga)\ : \ \ga\in\Ga\}$
in $\L (\ell^2(\Ga))$ for the operator norm.
Recall also that $\delta_e$ factorizes to a
faithful tracial state on $C^*_r(\Ga).$
The finite dimensional representations of 
$P\SL$ do not factorize through $C^*_r(P\SL)$,
since $P\SL$ is not amenable. As a consequence, 
Theorem~\ref{Theo1}
implies that $\delta_e$
is the unique  tracial state on $C^*_r(P\SL)$.
This also follows  from \cite{BCH},
where a different method is used.

 Theorem~\ref{Theo1} leaves open the
problem of existence of \emph{infinite}, semi-finite
traces on $C^*(SL_n(\ZZ)).$  We do not know
whether
such traces exist.
Using  \cite{BCH},
we can only show that no such trace
exists on $C_r^*(\PSL).$ 
In fact, this result is valid for a more 
general class of groups including $PSL_2(\ZZ)$
(see Proposition~\ref{Theo3} below).

This paper is organized as follows. Sections~\ref{S:Factors} and ~\ref{S:Subgroups}
are devoted to some general facts.
The proof of Theorem~\ref{Theo1} is spread over three sections:
in Section~\ref{S:Proof1}, we show that
the proof splits into two cases
which are then treated accordingly
in Sections~\ref{S:Proof2} and ~\ref{S:Proof3}. 
In Section~\ref{S:Superigidite}, we show that 
Theorem~\ref{Theo2} is a consequence of Theorem~\ref{Theo1}.
Corollary~\ref{Theo4} is proved in Section~\ref{S:Kirchberg}
and  Section~\ref{S:Theo3} is devoted to
a remark on the problem of the existence
of infinite traces. 

\par

\medskip\n
\textbf{Acknowlegments}
We are grateful to S. Popa who pointed out
to us  Connes' question  from  \cite{Jones}
and suggested to emphasize the superrigidity result
Theorem~\ref{Theo2}. Thanks are also due 
to E. Blanchard, E. Kirchberg, and P. de la Harpe,
 for interesting comments.

\section{Factor representations and characters}
\label{S:Factors}
We review some general facts
concerning  the relationships
between central  positive
definite  functions on groups and
factor representations.
Details can be found in  \cite[Chapters 6 and 17]{Dixmier}
or \cite{ThomaTyp1}.

Let $\Ga$ be a discrete group. We are interested
in representations of $\Ga$  in the unitary group
of a finite von Neumann algebra.

Recall that a finite
trace or a tracial state on a 
$C^*$--algebra $A$ with unit $1$
is a  linear functional
$\tau$ on $A$ which has the property
$$\tau (xy)= \tau (yx)\qquad\text{for all} \quad x, y\in A,$$
which is positive (that is, $\tau(x^*x)\geq 0$
for all $x\in A$),  and which is normalized by $\tau (1)=1$. The trace
$\tau$ is faithful if 
$\tau (x^*x)\neq 0$ for all $x\neq 0$.

Let $M$ be a finite von Neumann algebra, with faithful  trace normal $\tau.$
Let $\pi: \Ga\to U(M)$ be a group homomorphism.
The function $\vfi=\tau\circ \pi: \Ga\to \CCC$  
has the following properties:
\begin{itemize}
\item[(i)] $\f$ is central;
\item [(ii)] $\f$ is positive definite; 
\item [(iii)] $\f(e)=1.$
\end{itemize}

Let  $CP(\Ga)$ denote the set of functions $\f: \Ga\to \CCC$  
with Properties (i), (ii) and (iii) above. 

Ley  $\f\in CP(\Ga).$ Then there exist a finite von Neumann
algebra $M_\f,$ with faithful normal trace $\tau_\f,$ and 
a group homomorphism $\pi_\f: \Ga\to U(M_\f)$ such 
that $\f=\tau_\f\circ \pi_\f.$ Indeed, 
by GNS--construction, 
there exists a cyclic unitary representation $\pi_\f$ of $\Ga$ 
on a Hilbert space $\H_\f$ with a cyclic unit vector $\xi_\f$ 
such that
$$
\f (\g) = \langle\pi_\f (\g)\xi_\f, \ \xi_\f\rangle \qquad \text{for all}\quad  \ga\in \Ga.
$$
Since $\f$ is central, there exists another unitary representation
$\rho_\f$ of $\Ga$ on $\H_\f$ which commutes with $\pi_\f$ (that is,
$\pi_\f(\ga)\rho_\f(\ga')=\rho_\f(\ga')\pi_\f(\ga)$ for 
all $\ga,\ga'\in \Ga$) and with the property
that 
$$
\rho_\f(\ga)\xi_\f=\pi_\f(\ga^{-1})\xi_\f\qquad \text{for all}\quad \ga\in \Ga.
$$ 
Let $M_\f=\pi_\f(\Ga)''$ be the
von Neumann subalgebra of $\L(\H_\f)$ generated by $\pi_\f(\Ga),$
where 
$\CalS' = \{T\in \L(\H_\f)\ :\ TS=ST \tous S\in \CalS\}$ denotes
the commutant of a subset $\CalS$ of $\L(\H_\f).$
The mapping 
$$
T\mapsto \langle T\xi_\f, \ \xi_\f\rangle\tous T\in M_\f
$$ is a 
faithful normal trace $\tau_\f$
on $M_\f$ and $\f =\tau_\f\circ \pi_\f.$

Moreover, if $N_\f=\rho_\f(\Ga)''$ is the
von Neumann subalgebra of $\L(\H_\f)$ generated by $\rho_\f(\Ga),$ then 
$$
M_\f'=N_\f\et N_\f'=M_\f.
$$
 In particular, the common centre of $M_\f$ and $N_\f$ is $M_\f\cap N_\f.$

As an important example, let $\f=\delta_e$
be  the Dirac function 
 at  the group unit $e.$  Then $\f\in CP(\Ga).$ The  unitary representations $\pi_\f$ and $\rho_\f$
associated to $\f$ are  the left and
right regular representations $\lambda_\Ga$
and $\rho_\Ga$ on $\ell ^2(\Ga).$
Morever, $M_\f$ is the von Neumann algebra $L(\Ga)$
of $\Ga.$

The set $CP(\Ga)$ is a compact and convex subset of the vector space
of all bounded functions on $\Ga,$ equipped with the weak *-topology.
The set of extremal points $E(\Ga)$ of $CP(\Ga)$
is the set of all indecomposable central positive
definite  functions on $\Ga.$
By Choquet  theory, every
$\varphi\in CP(\Ga)$ may be written as a integral 
$$
\varphi=\int_{E(\Ga)}\psi d\mu (\psi)
$$
for a probability measure $\mu$ on $E(\Ga),$ at least
when $G$ is countable.
For $\f\in CP(\Ga),$ we have that $M_\f$ is a factor  if and only if 
$\vfi\in E(\Ga)$.
As an example, the Dirac function $\delta_e$ belongs to $E(\Ga)$ if and only 
if $\Ga$ is an ICC group.

Let $M$ be a finite von Neumann algebra,
with faithful normal trace $\tau,$ and let $\pi: \Ga\to U(M)$ be a homomophism
such that $\pi(\Ga)''=M.$ 
Observe that, if we set $\f=\tau\circ \pi\in CP(\Ga),$ then,  with the notation above,
the mapping $\pi_\f(\ga)\mapsto \pi (\ga)$ extends to an isomorphism
$M_\f\to M$ of von Neumann algebras.

A homomorphism
$\pi: \Ga\to U(M)$ for a finite factor $M$ such that $\pi(\Ga)''=M$
will be called a  finite factor representation of $\Ga.$ 
We say that two such representations $\pi_1:\Ga\to U(M_1)$ and
$\pi_2:\Ga\to U(M_2)$ are
quasi-equivalent if there exists an isomorphism 
$\Phi:M_1\to M_2$
such that $\Phi(\pi_1(\ga))=\pi_2(\ga)$ for all $\ga\in \Ga.$
Summarizing the discussion above, we see that $E(\Ga)$ classifies the finite
factor representations  of $\Ga$, up to quasi-equivalence.

 The set $E(\Ga)$ parametrizes also 
the indecomposable traces on the full $C^*$-algebra of 
$\Ga.$
 Recall that the full $C^*$-algebra $C^*(\Ga)$ of $\Ga$  is the 
 $C^*$-algebra with the universal
property that every unitary representation of $\Ga$ on a Hilbert 
space $\H$ extends to a $*$--homomorphism $C^*(\Ga)\to \L(\H).$
The algebra $C^*(\Ga)$ can be realized as
completion of the group algebra ${\CCC}[\Ga]$ under 
the norm 
$$
\left\Vert\sum_{\g\in\Ga}c_{\g}\g\right\Vert =\sup\left\{\left\Vert\sum_{\g\in\Ga}c_{\g}
\pi (\g)\right\Vert\ :\ \pi \in {\mathop{\rm Rep}}(\Ga)\right\},
$$
where ${\mathop{\rm Rep}}(\Ga)$ denotes the set of 
(equivalence classes of) cyclic unitary representations of $\Ga$.

We will view $\Ga$ as a subgroup
of the group of unitaries in $C^*(\Ga)$ by
means of the canonical embedding $\Ga\to C^*(\Ga).$
Every trace on $C^*(\Ga)$ defines by restriction
to $\Ga$ an element of $CP(\Ga).$ Conversely,
every $\f\in CP(\Ga)$ extends to a trace on $C^*(\Ga),$
since, as seen above, $\f(\ga)=  \langle\pi_\f (\ga)\xi_\f, \ \xi_\f\rangle$
and $\pi_\f$ is a unitary representation of $\Ga.$

\section{Some subgroups of $\SL$}
\label{S:Subgroups}
Let
$n$ is a fixed integer with $n\geq 2.$
For a pair of integers $(i,j)$ with $1\leq i\neq j\leq n,$ denote by $e_{ij}$ the 
corresponding elementary matrix, that is,
the $(n\times n)$-matrix  with $1$'s on the diagonal, $1$
at the $(i,j)$-entry,  and  $0$ elsewhere. 
It is well-known that $\SL$ 
is generated by 
$$\{e_{ij}\ :\ 1\leq i\neq j\leq n\}.$$
Moreover, for $n\geq 3,$ any two elementary matrices 
are conjugate
inside $\SL.$ Indeed,
observe that the matrix
 $$s_{ij}=e_{ij}e_{ji}^{-1}e_{ij}\in \SL$$
permutes the $i$-th and the $j$-th standard unit vectors
of $\ZZ^n,$ up to a sign. Hence, 
if $e_{kl}$ and $e_{pq}$ are two elementary matrices,
conjugation by a suitable product of  matrices of the form
$s_{ij}$ will carry  $e_{kl}$ into $e_{pq}$ or $e_{pq}^{-1}.$
Now, $e_{pq}$ and $e_{pq}^{-1}$
are conjugate  via a suitable diagonal matrix in $\SL,$
when $n\geq 3.$

The proof of the following two lemmas is by straightforward computation.
We will always view an element $a\in \ZZ^n$ as column vector.
Its transpose $a^t$ is then a row vector.
We denote by $e_1,\dots,e_n$ the standard unit vectors in $\ZZ^n.$

\begin{lemma}
\label{Lem-Centralisateur}
Let $k$ be a non-zero integer and let $i,j\in \{1,\dots,n\}$
with $1\leq i\neq j\leq n.$
The centralizer of $e_{ij}^k$  in $\SL$ consists of 
all matrices 
with  $\eps e_i$ as
$i$-th column and $\eps e_j^t$ as $j$-th row
for $\eps\in\{\pm 1\}.$ $\bsq$
\end{lemma}
For instance, the centralizer of $e_{12}^k$ is
the subgroup of all matrices of the form
$$
\begin{pmatrix}
\eps&*&*&\cdots&*\\
0&\eps&0&\ldots&0\\
0&*&*&\cdots&*\\
&\vdots&\vdots&\ddots&\\
0&*&*&*&*\\
\end{pmatrix},
$$
for $\eps\in\{\pm 1\}.$

For $j\in\{1,\dots, n\},$ let $V_j\cong \ZZ^{n-1}$ be the subgroup
generated by 
$$\{e_{ij}\ :\ 1\leq i\leq n, i\neq j\};$$
for instance, $V_1$ is the set of  matrices of the form
$$
\left(
\begin{array}{cccc} 
1&0&\ldots&0\\
*&1&\ldots&0\\
\vdots&&\ddots&\\
*&0&\ldots&1\\
\end{array}
\right).
$$
\begin{lemma}
\label{Lem-Normalisateur}
The normalizer of $V_j$ in $\SL$ is the subgroup
$G_j$ of all matrices in $\SL$
with  $\eps e_j$ as $j$-th row
for $\eps\in\{\pm 1\}.\bsq$
\end{lemma}
Thus, for instance, 
the normalizer $G_1$ of $V_1$ is the group of all matrices
$$
\begin{pmatrix}
\eps&0&\cdots&0\\
*&&&\\
\vdots&&A&\\
*&&&\\
\end{pmatrix},
$$
where $A\in GL_{n-1}(\ZZ)$ and $\eps=\det A.$

Up to a subgroup of index two,
$G_j$ is isomorphic to 
the semi-direct product
$SL_{n-1}(\ZZ)\ltimes\ZZ^{n-1}$ 
for the natural action of $SL_{n-1}(\ZZ)$ on $\ZZ^{n-1}.$

We will have also to consider the transpose subgroups
$V_i^t$ generated by  
$$\{e_{ij}\ :\  1\leq j\leq n, j\neq i\}.$$
Observe that $V_j\cap V_i^t$ is the copy of $\ZZ$ generated
by $e_{ij}$ for $i\neq j.$
The normalizer of  $V_i^t$ in $\SL$ is
of course the group $G_i^t.$
Observe also that 
$$V_j\subset G_i^t \qquad \text{and}\qquad V_i^t\subset G_j$$
for all $i\neq j.$

We will refer to subgroups of the form $V_j$ and $V_i^t$
as to the copies of $\ZZ^{n-1}$ inside $\SL.$

\section{Proof of Theorem \ref{Theo1}: A preliminary reduction}
\label{S:Proof1}
The starting point of the
proof of Theorem~\ref{Theo1} is the following classification
from  \cite[Proposition 9]{Burger}
of the measures on the $n$-dimensional torus
$\TT^n$ which  are invariant under the natural action
of $\SL$; for a more  elementary proof in the case
$n=2,$ see \cite{DaniKeane}.
\par
\begin{lemma}
 \textbf{([Bur])}
\label{Lem-Burger}
Let $n\geq 2$ be an integer.
Let $\mu$ be a $\SL$--invariant ergodic probability
measure on the Borel subsets of $\TT^n$. Then either $\mu$ is concentated on
a finite $\SL$--orbit or $\mu$ is the normalized Lebesgue measure 
on $\TT^n$.
\end{lemma}
Recall that a point $x\in {\TT}^n ={\RRR}^n/\ZZ^n$ has 
a finite $\SL$--orbit if and only if $x \in {\QQ}^n/\ZZ^n$.

Let $n\geq 3$ and let 
$$
\f: \SL\to \CCC
$$ be an indecomposable 
central positive definite 
function on 
$\SL$, fixed throughout the proof.

As in Section~\ref{S:Factors}, let $\pi$ and $\rho$ be the 
corresponding commuting factor representations
of $\Ga$ on the Hibert space $\H$
with cyclic vector $\xi$ such that
$$\f (\ga) = \langle\pi (\ga)\xi, \ \xi\rangle
= \langle\rho (\ga^{-1})\xi, \ \xi\rangle, \tous \ga\in \Ga.$$

Fix  any copy $V=V_j$ or $V=V_j^t$ of $\ZZ^{n-1}$ inside
$\SL$ and
consider the restriction  $\fV$ to $V$.

As $\f$ is central, $\fV$ is a
$G$--invariant positive 
definite function on $V$,
where 
$$G=G_j \qquad\text{or}\qquad G=G_j^t
$$ 
is the normalizer of $V$
in $\SL.$
 Since $G$ contains a copy of 
the semi-direct product $SL_{n-1}(\ZZ)\ltimes \ZZ^{n-1}$
(for the usual action in case $V=V_j$ and for the inverse transpose
of the usual action in case $V=V_j^t$), we have
$$
 \f (Ax) =\f (x)  
\qquad \text{for all} \quad x\in \ZZ^{n-1},\ \quad A\in \Sl.
$$
Thus, by Bochner's theorem, $\fV$ is the Fourier
transform of a $\Sl$--invariant probability 
measure on the torus 
$$\TT ^{n-1}\cong \widehat{V}.$$
Let $({\O}_i)_{i\geq 1}$ denote the sequence
of finite $\Sl$--orbits
in $\TT ^{n-1}$. For each $i\geq 1,$ denote by
$\mu_{{\O}_i}$ the uniform distribution
on $\O_i,$ that is, the 
 probability measure 
$$
\mu_{{\O}_i}=\frac{1}{|\O_i|} \sum_{\chi\in\O_i} \delta_{\chi}
$$
 on $\TT ^{n-1}.$ 
Lemma~\ref{Lem-Burger}  shows  that $\mu$ has a decomposition
as a convex combination
$$
\mu = t_\infty^{(V)} \mu_{\infty}+\sum_{i\geq 1}{ t_i^{(V)} \mu_{\O_i}}\qquad
\text{with}
\qquad t_\infty^{(V)}+\sum _{i\geq 1} t_i^{(V)}=1,\  t_\infty^{(V)}\geq 0, \ t_i^{(V)} \geq 0,
$$
where $\mu_{\infty}$ is the normalized Lebesgue measure
on $\TT ^{n-1}$.
Thus, we obtain a corresponding decomposition of $\fV$
$$
\fV = t_\infty^{(V)}\delta_e +\sum_{i\geq 1} t_i^{(V)} \psi_{\O_i}\qquad
\text{with}
\qquad t_\infty^{(V)}+\sum _{i\geq 1} t_i^{(V)}=1,\  t_\infty^{(V)}\geq 0, \ t_i^{(V)}\geq 0,
$$
where  $\psi_{\O_i}$ is the Fourier transform of
the measure $\mu_{\O_i}$.

By  general theory, we have
a corresponding decomposition 
of $\H$ into a direct sum of $\pi(V)$--invariant subspaces 
$$
\H =  \H_\infty^{V}\oplus \bigoplus_{\chi\in {\QQ}^{n-1}/\ZZ^{n-1}} \H_{\chi}^{V}
$$
where  $\H_\chi^{V}$ is the subspace  on which $V$  acts according to the character $\chi,$  that is,
$$
\H_\chi^{V}=\{\eta\in\H\: \ \pi(v)\eta=\chi(v)\eta \tous v\in V\}
$$
and where
$\H_\infty^{V}$ is a subspace 
on which $\pi(V)$ is a multiple of the regular representation
$\lambda_V$ of $V.$ 
Observe that some of these subspaces may be $\{0\}$.
Observe also that, since the representation $\rho$ commutes
with $\pi,$  each of the subspaces $\H_{\chi}^{V}$
and $\H_\infty^{V}$ is invariant under the whole
of $\rho(\SL).$

We claim that we have the following dichotomy.
\begin{lemma} 
\label{Lem-Dicho}
We have
\begin{itemize}
\item either $\H =   \bigoplus_{\chi\in {\QQ}^{n-1}/\ZZ^{n-1}} \H_{\chi}^{V}$
for every copy $V$ of $\ZZ^{n-1}$ in $\SL,$ or
\item $\H =  \H_\infty^{V}$
 for every copy $V$ of $\ZZ^{n-1}$ in $\SL$.
\end{itemize}
\end{lemma}
\begin{proof}
Let $V$ be a copy of $\ZZ^{n-1}$ with
$\H_\infty^{V}\neq \{0\}.$
We will show that 
$$
\H=\H_\infty^{W} \qquad\text{for every copy $W$ of $\ZZ^{n-1}$ in $\SL.$}
$$
Clearly, this will prove the lemma.

\n
 $\bullet$ {\it First step:} Let $W$ be a copy of $\ZZ^{n-1}$
for which we assume that
$V\cap W \neq \{0\}.$ We claim that 
$\H_{\infty}^{V}=\H_{\infty}^{W}.$

Indeed,  $V\cap W $
is the copy of $\ZZ$  generated 
by the appropriate elementary matrix.
We have two decompositions of $\H:$
$$
\H=\H_\infty^{V}\oplus \bigoplus_{\chi\in {\QQ}^{n-1}/\ZZ^{n-1}} \H_{\chi}^{V}
\qquad \text{and}\qquad 
\H=\H_\infty^{W}\oplus \bigoplus_{\chi\in {\QQ}^{n-1}/\ZZ^{n-1}} \H_{\chi}^{W}.
$$
Consider the restriction  of $\pi$ to 
$V\cap W.$ Each one of the subspaces   
$$
\bigoplus_{\chi\in {\QQ}^{n-1}/\ZZ^{n-1}} \H_{\chi}^{V}
\qquad \text{and} \qquad
\bigoplus_{\chi\in {\QQ}^{n-1}/\ZZ^{n-1}} \H_{\chi}^{W}
$$
has a decomposition into a direct sum
of subspaces under which  $\pi(V\cap W)$
acts according to a  character of $V\cap W.$

On the other hand, the representation $\pi|_{V\cap W}$ restricted to  
$\H_{\infty}^{V}$ or to $\H_{\infty}^{W}$
is a multiple of the regular representation $\lambda_{V\cap W},$
since $\lambda_{V}|_{V\cap W}$ and $\lambda_{w}|_{V\cap W}$
are mutiples of $\lambda_{V\cap W}.$
It follows that we  necessarily have $\H_{\infty}^{V}=\H_{\infty}^{W}.$

\n
 $\bullet$ {\it Second step:} Let $W$ be now an arbitrary copy of $\ZZ^{n-1}.$
We claim that we still have $\H_{\infty}^{V}=\H_{\infty}^{W}.$

Indeed, as is readily verified, we can find 
two copies $W^1$ 
and $W^2$ of $\ZZ^{n-1}$ with
$$
V\cap W^1 \neq \{0\},\quad W^1\cap W^2 \neq \{0\}, \quad\text{and}\quad
W^2\cap W \neq \{0\}.
$$
Therefore, by the first step, we have 
$$
\H_{\infty}^{V}=\H_{\infty}^{W^1}, \quad \H_{\infty}^{W^1}=\H_{\infty}^{W^2}
 \qquad \H_{\infty}^{W^2}= \H_{\infty}^{W},
$$
so that $\H_{\infty}^{V}=\H_{\infty}^{W}.$

\n
 $\bullet$ {\it Third step:} We claim that 
$\H_{\infty}^{V}=\H.$

Indeed,  by the second step,
we have  
$$\H_{\infty}^{V}=\H_{\infty}^{W}
\qquad\text{for every copy $W$ of $\ZZ^{n-1}$ in $\SL.$} .$$ 
Since $\H_{\infty}^{W}$ is invariant under
$\pi(W),$ it follows that $\H_{\infty}^{V}$
is invariant under
$\pi(\SL).$ 

On the other hand, $\H_{\infty}^{V}$ is also invariant under
$\rho(\SL).$ Since $\pi$ is a factor representation
and since  $\H_{\infty}^{V}\neq \{0\},$
the claim follows.
$\bsq$
\end{proof}

We have to consider separately the two possible
decompositions of $\H$ given by the previous lemma.
We will see 
that the first one corresponds to 
a character of a congruence quotient,
and that the second one to a character
induced from the centre.

\section{Proof of Theorem \ref{Theo1}: First case}
\label{S:Proof2}
With the notation from the last section,
we assume in this section that
 $$
\H =   \bigoplus_{\chi\in {\QQ}^{n-1}/\ZZ^{n-1}} \H_{\chi}^{V}
\qquad\text{for every copy $V$ of $\ZZ^{n-1}$ in $\SL.$}
$$
We claim that there exists some integer $N\geq 1$ such that
$\pi$ is trivial  on the congruence normal  subgroup
$$
\Ga(N) = \{\ga\in\SL\ :\ \ga \ \equiv I \mod N\}.
$$

Let $\ga_0,\ga_1,\dots,\ga_{d}$ denote the 
elementary matrices in $\SL,$
where $d=n(n-1)-1.$

For every $k\in \{0,\dots, d\},$
 we  have a decomposition 
$$
\H =   \bigoplus_{\alpha\in \QQ/\ZZ}  \H_{\alpha}^{\ga_k}
$$
of $\H$ under the action
of the unitary operator $\pi(\ga_k),$ where 
$\H_{\alpha}^{\ga_k}$
is the eigenspace (possibly equal to $\{0\}$) of $\pi(\ga_k)$ corresponding to $\alpha.$

\begin{lemma}
\label{Lem-FiniteDec}
There exists an integer $N\geq 1$
such that  $\pi(\ga_0^N),\pi(\ga_1^N),\dots,\pi(\ga_d^N) $ 
have a non-zero common invariant vector in $\H.$
\end{lemma}
\begin{proof}
Let $M$ be the factor
generated by $\pi(\Ga)$ and denote by $\tau$
the trace on $M$ defined by $\f.$

Write the elements in $\QQ/\ZZ $
as a sequence $\{\alpha_i)_{i\geq 1}.$ 
For every $i\geq 1,$ let 
$$p_i:\H\to \H_{\alpha_i}^{\ga_0}$$
 denote the orthogonal projection
onto  $\H_{\alpha_i}^{\ga_0}.$
Observe that $p_i\in M$ (in fact,
$p_i$ belongs to the abelian von Neumann
algebra generated by $\pi(\ga_0)$).
We have $\tau(p_i)\in [0,1]$
and
 $\sum_{i\geq 1}\tau (p_i)=1,$
since $\sum_{i\geq 1} p_i= I,$

Let  $\eps$ be a real number 
with  
$$0<\eps<1/2^d.$$
There exists $i_0\geq 1$ such that  
$$
\sum_{i=1}^{i_0}\tau(p_i)\geq 1-\eps.
$$
Since elements in $\QQ/\ZZ$ have finite
order, we can find an integer $N\geq 1$ such  that
$$\alpha_i^N=1 \tous i\in \{1,\dots, i_0\}.
$$ 
Then $\pi(\ga_0^N)$ acts as the identity
on  
$$\bigoplus_{i=1}^{i_0}  \H_{\alpha_i}^{\ga_0}.$$

For  $l\in \{0,1,\dots, d\},$
let $\H^{\ga_l^N}$ be the subspace of 
$\pi(\ga_l^N)$--invariant vectors
in $\H.$ We claim that 
$$
\H^{\ga_0^N}\cap \H^{\ga_1^N}\cap\cdots\cap  \H^{\ga_d^N}\neq\{0\}. 
$$ 

For every $k\in \{0,1,\dots, d\},$ let 
$q_k$ denote the orthogonal projection
onto  
$$
\H^{\ga_0^N}\cap \H^{\ga_1^N}\cap\cdots\cap \H^{\ga_k^N}.
$$
It is clear that $q_k\in M.$ 

We claim that 
$$
 \tau (q_k) \geq 1-2^k \eps \qquad \text{for all}\quad k=0,1,\dots, d.\leqno{(1)}
$$ 
Once proved, this will imply that
$$\tau (q_d)\geq 1-2^d\eps >0,$$
and hence $q_d\neq 0$ since $\tau$ is faithful on $M;$
this will finish the proof
of the lemma. 

To prove $(1),$ we proceed by induction on
$k.$ Since 
$$
\bigoplus_{i=1}^{i_0}  \H_{\alpha_i}^{\ga_0}\subset \H^{\ga_0^N},
$$
we have $q_0\geq \sum_{i=1}^{i_0}p _i.$ Hence,
$$
\tau (q_0)\geq \sum_{i=1}^{i_0}\tau(p _i)\geq 1-\eps,
$$
and this proves $(1)$ in the case $k=0.$

Let $k\geq 1$ and assume that
$$
\tau (q_{k-1}) \geq 1-2^{k-1}\eps.\leqno{(2)}
$$ 
Set
$$
\K=\H^{\ga_0^N}\cap \H^{\ga_1^N}\cap\cdots\cap \H^{\ga_{k-1}^N}
$$
and set $q= q_{k-1},$  the orthogonal projection
on $\K$.
 
Since any two elementary matrices
are conjugate, we have $\ga_k= s\ga_0 s^{-1}$ for
some element $s\in \SL.$ 
Observe that 
$$
\H^{\ga_k^N}= \pi(s)\H^{\ga_0^N}.
$$
Consider the operator  
$$T= (1-q)\pi(s^{-1}) q$$ 
on $\H.$
Observe that $T\in M.$ For $\eta\in \H,$
we have $T(\eta) =0$ if and only if 
$\pi(s^{-1}) q(\eta)\in \K,$
that is, if and only $q(\eta)\in \pi(s)\K.$
Hence
$$
  \Ker T= (\K \cap \pi(s)\K)\oplus \K^\perp.\leqno{(3)}
$$
Let 
$$p_{\Ker T}:\H\to \Ker T $$ 
be the orthogonal projection 
on $\Ker T.$ Then $p_{\Ker T}\in M,$ since $T\in M.$
Moreover, since the range of $T$ is contained
in $\K^{\perp},$ we have
$$
\tau(1-q) \geq \tau(I)-\tau (p_{\Ker T}) =1-\tau (p_{\Ker T}).
$$
Hence, by $(2)$, 
$$
\tau (p_{\Ker T})\geq 1-2^{k-1}\eps.\leqno{(4)}
$$
We have, by $(3)$
$$
\tau (p_{\Ker T}) =\tau(p_{\K \cap \pi(s)\K}) + \tau (1-q),
$$
where $p_{\K \cap \pi(s)\K}\in M$ is the orthogonal projection
on $\K \cap \pi(s)\K.$
Now,
$$
\K \cap \pi(s)\K\subset \K \cap \pi(s)\H^{\ga_0^N} =\K\cap \H^{\ga_k^N}.
$$
Since $q_k$ is the orthogonal projection
on $\K\cap\H^{\ga_k^N},$ it follows in view
of $(2)$ and $(4)$ that
\begin{eqnarray*}
\tau (q_k)&\geq& \tau(p_{\K \cap \pi(s)\K}) \\
 &=&\tau (p_{\Ker T})-(1-\tau(q))\\
&\geq& (1-2^{k-1}\eps)-2^{k-1}\eps=1-2^{k}\eps.
\end{eqnarray*}
This proves the claim $(1)$ and finishes the proof of the lemma.
 $\bsq$
\end{proof}
\begin{corollary}
\label{Cor-FiniteDec}
Under the assumption made at the beginning of
this section,
there exists an irreducible  representation
$\pi_0$ of the congruence quotient
$$\SL/ \Ga(N^2)\cong SL_n(\ZZ/N^2\ZZ)$$ such that  
$\f$ is the (normalized) character of $\pi_0$ lifted to $\SL,$
 where $N$ is as in Lemma~\ref{Lem-FiniteDec}.
\end{corollary}
\begin{proof}
By the previous lemma, the subspace $\K$
of the common  invariant vectors  under 
$\pi(\ga_0^N),\pi(\ga_1^N),\dots,\pi(\ga_d^N) $ 
is non-zero. Let $\Ga$ be the subgroup
of $\SL$ generated by 
$$\{\ga_0^N,\ga_1^N,\dots,\ga_d^N\}.$$
By \cite[Proposition~2]{Tits}, $\Ga$ 
contains the congruence normal subgroup
$\Ga(N^2).$ 

Consider the subspace
$$
\H^{\Ga(N^2)}=\{\eta\in\H\ :\ \pi(\ga)\eta=\eta\tous \ga\in\Ga(N^2) \},
$$ 
of 
$\pi(\Ga(N^2))$-invariant vectors.
Then $\H^{\Ga(N^2)}\neq \{0\}$ since $\K\subset \H^{\Ga(N^2)}.$
Moreover, 
$\H^{\Ga(N^2)}$ is invariant under  
 $\pi(\SL),$ as $\Ga(N^2)$ is normal in $\SL.$ 

On the other hand,  $\H^{\Ga(N^2)}$ 
is clearly invariant under $\rho(\SL).$
It follows that 
$$\H^{\Ga(N^2)}=\H.$$
Hence, $\pi$ factorizes
through the finite group
$\SL/ \Ga(N^2).$
It follows that $\H$ is finite-dimensional,
that $\pi$ is a
equivalent to a multiple $m\pi_0$ of an irreducible
representation $\pi_0$ of $\SL/ \Ga(N^2)$
with $m=\dim (\pi_0),$ and that $\f$
is the normalized character of $\pi_0.$
  $\bsq$
\end{proof}

\section{Proof of Theorem \ref{Theo1}: Second case}
\label{S:Proof3}
With the notation as in Section~\ref{S:Proof1},
we assume now that
 $$
\H =  \H_{\infty}^{V}
\qquad\text{for every copy $V$ of $\ZZ^{n-1}$ in $\SL.$}
$$ 
This is equivalent to:
$$\f|_V= \delta_e 
\qquad\text{for every copy $V$ of $\ZZ^{n-1}$ in $\SL.$}
$$

Let $\chi_{\f}$ be the unitary character
of the centre $C=\{\pm I\}$ of $\SL$ such that
$$
\f(z\ga)=\chi_{\f}(z)\f(\ga) \text{for all}\quad z\in C,\ga\in \SL.
$$
We claim that 
$$
\f(\ga)=
\begin{cases}
0&\text{if $\ga\in \SL\setminus C$}\\
\chi_{\f}(\ga) &\text{if $\ga\in C$.} 
\end{cases}
$$

The following proposition, which is  of independent
interest, will play a crucial r\^ole.
\begin{proposition}
\label{Pro-Decomposition}
Every matrix  $\ga\in \SL$ 
 is conjugate to the product $g_1g_2g_3$ 
of three matrices of the form
$$
g_1= \begin{pmatrix} 
1&*&\cdots&*\\
0&*&\cdots&*\\
\vdots&&*&\\
0&*&\cdots&*\\
\end{pmatrix} \in G_1^t\ , 
\qquad 
g_2= \begin{pmatrix} 
*&*&\cdots&*\\
\vdots&&*&\\
*&*&\cdots&*\\
0&0&0&1\\
\end{pmatrix}\in G_n \
$$
and 
$$
g_3= \begin{pmatrix} 
1&0&0&\cdots&0\\
*&1&0&\cdots&0\\
\vdots&\vdots&\vdots&\vdots&\vdots\\
*&0&0&\cdots&1\\
\end{pmatrix}\in V_1
$$
\end{proposition}
\begin{proof}

\n
 $\bullet$ {\it First step:} 
We first claim that $\ga$ is conjugate to a matrix
$\ga_1$ with first column of the form
$(*,0,*,0,\ldots,0)^t.$ 
This is Lemma~1 in \cite{Brenner}.
The result is proved by conjugating $\ga$ by
 permutation matrices (with sign ajusted) and by 
elementary matrices of the type $e_{ij}$
with $1<i\neq j\leq n.$ 

So, we can assume that the first column of $\ga$ is 
of the form $(k,0, l,0,\ldots,0)^t$ for $k,l\in\ZZ.$

\medskip
\noindent
$\bullet$ {\it Second step:} There exists a matrix
$\ga_1\in G_n$ such that the first column
of $\ga_1 \ga$ is  $(k,1,l ,0,\ldots,0)^t.$
Indeed, since $\gcd (k, l)=1,$ there exist
$p,q\in \ZZ$ such that $pk+ ql=1.$
We can take
$$\ga_1 =
\begin{pmatrix} 
1&0&0&\cdots&0\\
p&1&q&\cdots&0\\
\vdots&\vdots&\vdots&\vdots&\vdots\\
0&0&0&\cdots&1\\
\end{pmatrix}\in G_n.
$$

\medskip
\noindent
$\bullet$ {\it Third step:} There exists a matrix
$\ga_2\in G_1^t\cap G_n$ such that the first column
of $\ga_2\ga_1 \ga$ is  $(1,1, l,0,\ldots,0)^t.$
Indeed, we can take 
$$
\ga_2=
\begin{pmatrix} 
1&1-k&0&\cdots&0\\
0&1&0&\cdots&0\\
\vdots&\vdots&\vdots&\vdots&\vdots\\
0&0&0&\cdots&1\\
\end{pmatrix}\in  G_n.
$$

\medskip
\noindent
$\bullet$ {\it Fourth step:} There exists a matrix
$\ga_3\in V_1$ such that the first column
of $\ga_3\ga_2\ga_1 \ga$ is  $(1,0, 0,0,\ldots,0)^t.$
Indeed, we can take 
$$\ga_3=\begin{pmatrix} 
1&0&0&\cdots&0\\
-1&1&0&\cdots&0\\
-l&0&1&\cdots&0\\
\vdots&\vdots&\vdots&\vdots&\vdots\\
0&0&0&\cdots&1\\
\end{pmatrix}\in V_1.
$$

 By the last
step, $\ga_4=\ga_3\ga_2\ga_1 \ga\in G_1^t.$
We have 
$$ 
\ga_4\ga \ga_4^{-1}= \ga_4 (\ga_1^{-1}\ga_2^{-1})\ga_3^{-1}.
$$
The claim follows, since $\ga_1^{-1}\ga_2^{-1}\in G_n$
and $\ga_3^{-1}\in V_1.$
$\bsq$
\end{proof}
\begin{remark}
\label{Rem-Conjugate}
In the case $n\geq 4,$ 
 the previous proposition
can be improved:
every $\ga\in\SL$ is conjugate
to a product $g_1g_2\in G_1^tG_n.$
Indeed, in this case, 
the matrix $\ga_3$ in the fourth step of the proof
belongs to $G_n$ and hence 
$$
\ga_4\ga \ga_4^{-1}= \ga_4 (\ga_1^{-1}\ga_2^{-1}\ga_3^{-1})\in G_1^tG_n.
$$
\end{remark}

\bigskip
Returning to the proof of Theorem~\ref{Theo1},
 the previous proposition implies that 
it suffices to show that
$$
\vfi(\ga)=0\tous \ga\in G_1^tG_n V_1 \quad\text{with}\quad \ga\notin C.
$$
For this, several preliminary steps will
be needed.

We will use several times the following elementary lemma.
\begin{lemma}
\label{Lem-Hilbert}
 Let $\Ga$ be a group and
$\Hpi$ a unitary representation of $\Ga.$ 
Let $\psi=\langle \pi(\cdot)\xi,\xi\rangle$
be an associated positive definite function
such that $\psi =\delta_e.$ Then, for every 
sequence $(g_k)_{k\in\NN}$ of pairwise distinct
elements $g_k\in \Ga,$ the sequence
$(\pi(g_k)\xi)_{k\in\NN}$ converges weakly
to $0$ in $\H.$
\end{lemma}
\begin{proof}
 For $k,l\in \NN$ with $k\neq l,$ we have
\begin{eqnarray*}
\langle \pi(g_k)\xi,\pi(g_l)\xi\rangle
&=&\langle \pi(g_l^{-1}g_k)\xi,\xi\rangle\\
&=&\psi(g_l^{-1}g_k) =0.
\end{eqnarray*}
Therefore, $(\pi(g_k)\xi)_{k\in\NN}$ 
is an orthonormal sequence in $\H$
and the claim follows.$\bsq$
\end{proof}

The first step in this part of the proof
of Theorem~\ref{Theo1} is to show that
$$
\vfi(\ga)=0\tous \ga\in G_1^t\cup G_n  \quad\text{with}\quad \ga\notin C.
$$

For elements $x,y$ in a group, let $[x,y]$
denote the commutator $x^{-1} y^{-1}xy.$
\begin{lemma}
\label{Lem-ExtNorm}
Let $V$ be a copy of  $\ZZ^{n-1}$ in $\SL$
and let $G$ be the normalizer
of $V.$ Then
$\f(\ga) = 0$ for every $\ga\in G\setminus C.$
\end{lemma}
\begin{proof}
  Write 
$$V=\ZZ x_1\oplus \cdots \oplus\ZZ x_n,$$
where $x_1,\dots,x_n$ are the elementary matrices
contained in $V.$

Let $\ga\in G\setminus C.$ We claim that there exists $i\in \{1,\dots, n\}$
such that
$$
 x_i^{-k}\ga x_i^k\neq x_i^{-l}\ga x_i^l\qquad \text{for all}
\quad k,l\in \ZZ,\ k\neq l.
$$ 
Indeed, otherwise there would exist
non-zero integers $k_i$ such
that $\ga$ is in the  centralizer
of  $x_i^{k_i}$ for all $i\in \{1,\dots, n\}.$
This would imply that $\ga\in C$
(see Lemma~\ref{Lem-Centralisateur}).

The  commutators $[\ga,x_i^k]$
belong to $V$ and are pairwise
distinct.
Hence, by Lemma~\ref{Lem-Hilbert}, the sequence
$(\pi ([\ga,x_i^k])\xi)_{k\in \NN}$
is weakly convergent to $0$ in $\H.$ 
For $k\in \NN,$  we have
\begin{eqnarray*}
\f(\ga)&=&\f(x_i^{-k}\ga x_i^k)\\
&=&\f(\ga [\ga,x_k])\\
&=&\langle \pi([\ga,x_i^k])\xi,\pi(\ga^{-1})\xi\rangle.
\end{eqnarray*}
Hence,
$$
\f(\ga)=\lim_k \langle \pi([\ga,x_i^k])\xi,\pi(\ga^{-1})\xi\rangle=0,
$$
as claimed.
$\bsq$
\end{proof}

The next step is to show that 
$$
\vfi(\ga)=0\tous \ga\in G_1^t G_n \quad\text{with}\quad \ga\notin C.
$$

\begin{lemma}
\label{Lem-ExtProduit}
Let $V, W$ be two copies of  $\ZZ^{n-1}$ in $\SL$
with $V\cap W\neq \{0\}.$
Let $G, H$ be the normalizers
of $V$ and $W,$ respectively. 
Let $\ga=gh$ with $g\in G,$
$h\in H,$ and $\ga\notin C.$ Then
$\f(\ga) = 0.$ 
\end{lemma}
\begin{proof}
If $g\in C$ or $h\in C,$
then $\ga\in G$ or $\ga\in H$
and then  $\vfi(\ga)=0,$
by Lemma~\ref{Lem-ExtNorm}.
Hence, we can assume that $g\notin C$
and $h\notin C.$

Let $x$ denote the elementary matrix 
such that  
$$
V\cap W =\langle x\rangle.
$$
It is readily verified that, for $k\in \ZZ\setminus\{0\},$ the centralizer of 
$x^k$ is contained in $G\cap H.$
Hence, we can assume that
$\ga$ does not belong to this centralizer,
that is, that  the elements $x^{-k}\ga x^k$
are pairwise distinct.

We have 
\begin{eqnarray*}
x^{-k}\ga x^k =x^{-k}g x^k x^{-k}h x^k
 =g[g,x^k]x^{-k}h x^k,
\end{eqnarray*}
Set $y_k=[g,x^k]x^{-k}h x^k.$
Observe that 
$V\subset G\cap H.$
Since $[g,x^k]\in V,$  we have $y_k\in H.$
Moreover, the elements $y_k$ are pairwise distinct,
since 
$$
y_k=g^{-1}x^{-k}\ga x^k.
$$
Hence, again by Lemma~\ref{Lem-Hilbert}, 
the sequence
$(\pi (y_k)\xi)_{k\in \NN}$
is weakly convergent to $0$ in $\H.$ 
As in the previous lemma, it follows that
$$
\f(\ga)=\lim_k \f(x^{-k}\ga x^k)=\lim_k \langle \pi(y_k)\xi,\pi(g^{-1})\xi\rangle=0. \ \bsq
$$
\end{proof}

We will also need the following consequence of Lemma~\ref{Lem-ExtProduit}.
\begin{lemma}
\label{Lem-ExtProduit2}
Let $V, W$ two copies of  $\ZZ^{n-1}$ in $\SL$
with $V\cap W\neq \{0\}.$
Let $G, H$ be the normalizers
of $V$ and $W,$ respectively. 
Let $(\ga_k)_{k\in \NN}$
be a sequence of  pairwise distinct elements in  $GH.$ 
Then
$(\pi(\ga_k)\xi)_{k\in \NN}$ converges weakly to $0$ in $\H.$
\end{lemma}
\begin{proof}
Observe that $(\pi(\ga_k)\xi)_{k\in \NN}$
is a bounded sequence in $\H.$ 
Therefore, it suffices to show that 
every subsequence  $(\pi(\ga_{k_i})\xi)_{i\in \NN}$ 
of $(\pi(\ga_k)\xi)_{k\in \NN}$
has a subsequence which weakly converges to $0.$

For $i\in \NN,$ write 
$\ga_{k_i}= g_{k_i} h_{k_i}$ for 
$g_{k_i}\in G$ and $h_{k_i}\in H.$

Since $C$ is finite and since the elements $\ga_{k_i}$
are pairwise distinct,
we can find a subsequence of $(\ga_{k_i})_i$, still denoted
by   $(\ga_{k_i})_i,$ 
such that 
$\ga_{k_j}^{-1}\ga_{k_i}\notin C$ for all $i\neq j.$
It follows that
$$
g_{k_j}^{-1}g_{k_i} h_{k_i}h_{k_j}^{-1}\notin C \tous i\neq j.
$$
From   Lemma~\ref{Lem-ExtProduit}, we deduce that, for all 
$i\neq j,$ 
\begin{eqnarray*}
\f(\ga_{k_j}^{-1}\ga_{k_i})&=& 
\f(h_{k_j}^{-1}g_{k_j}^{-1}g_{k_i} h_{k_i})\\
&=&\f(g_{k_j}^{-1}g_{k_i} h_{k_i}h_{k_j}^{-1})\\
&=& 0,
\end{eqnarray*}
since $g_{k_j}^{-1}g_{k_i}\in G$ and $ h_{k_i}h_{k_j}^{-1}\in H.$
As in the proof of Lemma~\ref{Lem-Hilbert},
this shows that 
$(\pi(\ga_{k_i})\xi)_{i}$  weakly converges to $0.$
$\bsq$
\end{proof}

We can now conclude the proof
of Theorem~\ref{Theo1}.
Let $\ga\in\SL\setminus C.$ We want to show that
$\f(\ga)=0.$

By Proposition~\ref{Pro-Decomposition}, we can 
assume that $\ga=g_1g_2g_3$ 
for matrices of the form
$$
g_1= \begin{pmatrix} 
1&*&\cdots&*\\
0&*&\cdots&*\\
\vdots&&*&\\
0&*&\cdots&*\\
\end{pmatrix} \in G_1^t\ , 
\qquad 
g_2= \begin{pmatrix} 
*&*&\cdots&*\\
\vdots&&*&\\
*&*&\cdots&*\\
0&0&0&1\\
\end{pmatrix}\in G_n \
$$
and 
$$
g_3= \begin{pmatrix} 
1&0&0&\cdots&0\\
a_2&1&0&\cdots&0\\
\vdots&\vdots&\vdots&\vdots&\vdots\\
a_n&0&0&\cdots&1\\
\end{pmatrix}\in V_1.
$$

If $g_3\in G_n,$ then $\ga$
is a non-central element in $G_1^tG_n,$
and it follows from Lemma~\ref{Lem-ExtProduit}
that $\f(\ga)=0.$
 We can therefore assume that $g_3\notin G_n,$
that is, $a_n\neq 0.$

Let $x$ be the elementary matrix $e_{2,n}$,
thus
$$
x= \begin{pmatrix} 
1&0&0&\cdots&0\\
0&1&0&\cdots&1\\
\vdots&\vdots&\vdots&\vdots&\vdots\\
0&0&0&\cdots&1\\
\end{pmatrix}.
$$
Then $x\in G_1^t\cap G_n$ and the
centralizer of every power $x^k$ 
for $k\neq 0$
is contained in $G_n.$ 
Hence, if $\ga$ is contained in the centralizer
of some power $x^k$ for $k\neq 0$, the claim  follows
from Lemma~\ref{Lem-ExtNorm}.
We can therefore assume that 
$$
x^{-k}\ga x^k\neq x^{-l}\ga x^l\qquad\text{for all}\quad k\neq l.
$$
We compute that
$$
x^{-k}g_3 x^k =
\begin{pmatrix} 
1&0&0&\cdots&0\\
a_2+ka_n&1&0&\cdots&0\\
a_{3}&0&1&\cdots&0\\
\vdots&\vdots&\vdots&\vdots&\vdots\\
a_n&0&0&\cdots&1\\
\end{pmatrix}.
$$
Hence $x^{-k}g_3 x^k=\alpha_k \beta,$ where
$$
\alpha_k= 
\begin{pmatrix} 
1&0&0&\cdots&0\\
ka_n&1&0&\cdots&0\\
0&0&1&\cdots&0\\
\vdots&\vdots&\vdots&\vdots&\vdots\\
0&0&0&\cdots&1\\
\end{pmatrix} \qquad\text{and}\qquad 
\beta= 
\begin{pmatrix} 
1&0&0&\cdots&0\\
a_2&1&0&\cdots&0\\
a_3&0&1&\cdots&0\\
\vdots&\vdots&\vdots&\vdots&\vdots\\
a_n&0&0&\cdots&1\\
\end{pmatrix}.
$$
Observe that $\alpha_k\in G_n$ for every $k.$
We have
\begin{eqnarray*}
x^{-k}\ga x^k&=&x^{-k}g_1g_2g_3 x^k\\
&=&(x^{-k}g_1x^k) (x^{-k}g_2x^k) (x^{-k}g_3 x^k)\\
&=&(x^{-k}g_1x^k) (x^{-k}g_2x^k)\alpha_k \beta.
\end{eqnarray*}
Now, since $x\in G_1^t\cap G_n,$
we have  $x^{-k}g_1x^k\in G_1^t$ and $x^{-k}g_2x^k\alpha_k\in G_n.$
It follows that 
$$
x^{-k}\ga x^k\beta^{-1} \in G_1^t G_n \qquad\text{for every}\quad k.
$$
Set 
$$\ga_k= x^{-k}\ga x^k\beta^{-1}.$$
Since $\ga$ is not in the centralizer
of $x^{k},$ we have $\ga_k\neq \ga_l$ for all $k\neq l.$
Hence, by Lemma~\ref{Lem-ExtProduit2}, the
sequence $(\pi(\ga_k))_{k\in \NN}$
converges weakly to $0.$
It follows that
\begin{eqnarray*}
\f (\ga)&=& \lim_k\f(\beta x^{-k}\ga x^k \beta^{-1})\\
&=& \lim_k\f(\beta\ga_k)\\
 &=& \lim_k\langle \pi(\beta\ga_k)\xi, \xi\rangle\\
&=& \lim_k\langle \pi(\ga_k)\xi, \pi(\beta^{-1})\xi\rangle\\
&=& 0.
\end{eqnarray*}
This concludes the proof of Theorem~\ref{Theo1}.$\bsq$

\section{Deducing Theorem~\ref{Theo2} from Theorem~\ref{Theo1}}
\label{S:Superigidite}
Let $\Ga=SL_n(\ZZ)$ for $n\geq 3.$
Let $M$ be a  finite factor, with trace $\tau,$ 
and let 
$\pi:\Ga\to U(M)$ be a group homomorphism
such that $\pi(\Ga)''=M.$
Then $\vfi=\tau\circ\pi$ is a character of $\Ga.$

Assume that  $M$ is finite dimensional.
Let $\pi_\f:\Ga\to U(M_\f)$ be the finite
factor representation associated to $\f$
(see Section~ \ref{S:Factors}).
The mapping $\pi_{\f}(\ga) \mapsto \pi(\ga)$
extends to an isomorphism $M_\f\to M$ of von Neumann algebras.
Hence $M_\f$ is finite dimensional and,
by  Theorem~~\ref{Theo1}, $\f$ is the character
 of an irreducible  finite dimensional
 representation of some congruence quotient
$SL_n({\ZZ}/N{\ZZ})$ for $N\geq 1.$
It follows that $\pi$ factorizes through $SL_n({\ZZ}/N{\ZZ}).$

Assume now that $M$ is infinite dimensional. 
By Theorem~~\ref{Theo1}, we have $\vfi=\widetilde\chi$ 
for a character $\chi$ of the centre $C.$ 
If $n$ is odd,  let $\Lambda=\Ga$ and, if $n$ is even,
let $\Lambda=\Ga(N)$ be a congruence subgroup for
$N\geq 3.$
Then $\Lambda$ has finite
index in $\Ga$ and $\Lambda\cap C=\{e\}.$
We  therefore  have $\vfi|_\Lambda=\delta_e.$
The GNS-representation of $\Lambda$
corresponding to $\delta_e$ is the
 regular representation $\lambda_{\Lambda}$
which generates the von Neumann algebra $L(\Lambda).$ 
The mapping $\lambda_{\Lambda}(\ga) \mapsto \pi(\ga)$
extends to a normal homomorphism $L(\Lambda)\to M.$

\begin{remark}
\label{Rem-Extension}
Observe that the conclusion in (ii) of Theorem~\ref{Theo2} is that
$\pi|_\Lambda$ extends to $L(\Lambda)$ and not just to $U(L(\Lambda)).$
P.~de la Harpe pointed out to me 
that this is a stronger statement: a homomorphism $U(M_1)\to U(M_2)$ between 
the unitary groups of two  finite factors $M_1,M_2$ does not necessarily
extend to an algebra homomorphism $M_1\to M_2.$
As a simple example, take $M_1=M_2(\CCC)$ and $M_2= M_4(\CCC)\cong M_2(\CCC)\otimes M_2(\CCC).$
The group homomorphism  $\pi: U(2)\to U(4), g\mapsto g\otimes g$ does not extend
to an algebra homomorphism $M_2(\CCC)\to M_4(\CCC).$
\end{remark}

\section{A question of Kirchberg}
\label{S:Kirchberg}

A conjecture of  Kirchberg  \cite[Section 8, (B4)]{Kirchberg}
is:

\medskip
The full $C^*$-algebra $C^*(SL_2(\ZZ)\times SL_2(\ZZ))$
of the direct product $SL_2(\ZZ)\times SL_2(\ZZ)$ has a faithful
tracial state.

\medskip
\n
As shown in \cite{Kirchberg},
 this problem is in fact equivalent to
 a series of outstanding  conjectures,
among them the following one which was suggested
by Connes in \cite[page 105]{Connes}:

\medskip
 Every  factor  of type $II_1$ with separable predual
is a subfactor of the ultrapower $R_{\omega}$  of the hyperfinite factor $R$ of
type  $II_1$. 
\medskip

\n
 A positive answer to the following
question of Kirchberg   \cite[Remark~8.2]{Kirchberg}
would imply 
the conjecture above:

\medskip
Does $C^*(SL_4(\ZZ))$ have a faithful tracial state?
\medskip

\n
Indeed, $SL_2(\ZZ)\times SL_2(\ZZ)$ embedds as a subgroup of 
$SL_4(\ZZ),$ for instance, via the mapping
$$
SL_2(\ZZ)\times SL_2(\ZZ)\ni(\ga_1,\ga_2)\to 
\begin{pmatrix}
\ga_1&0\\
0&\ga_2
\end{pmatrix}\in SL_4(\ZZ).
$$
 
A faithful tracial state on  $C^*(SL_4(\ZZ))$
would give, by restriction, a faithful tracial state on  
$C^*(SL_2(\ZZ)\times SL_2(\ZZ)).$

We proceed to show that
the answer to this question is negative.
In fact, the following stronger
result will be proved.
 We will consider the  copy 
$$
\Lambda=\left\{
\begin{pmatrix}
\ga&0\\
0&I
\end{pmatrix}
\ :\ \ga\in SL_2(\ZZ)\right \}\cong SL_2(\ZZ)
$$
of $SL_2(\ZZ)$ inside  $SL_n(\ZZ).$
\begin{corollary}
\label{Cor-Kirch}
Let $n\geq 3$ and set $\Ga=SL_n(\ZZ).$
Let $\f$ be a tracial state  on $C^*(\Ga).$
Then $\f|_{C^*(\Lambda)}$  is not faithful.
\end{corollary}
\begin{proof}
Let $\pi$ be the 
cyclic unitary representation of $\Ga$ corresponding to $\f.$
By Theorem~\ref{Theo1}, $\pi$ decomposes
as a direct sum 
$$
\pi_{\infty}\oplus\bigoplus_i \sigma_{i}\,,
$$
where
$\pi_{\infty}$ is a multiple of the regular representation
$\lambda_{\Ga},$ and 
where every representation $\sigma_{i}$ 
factorizes through
some congruence quotient $\Ga/\Ga(N_i).$ 

Let ${\rm Rep}_{\rm cong}(\Ga)$ denote
the set of all unitary representations of $\Ga$
which factorize through
some congruence quotient. 
In fact, as a consequence
of the positive answer to the congruence subgroup
problem, ${\rm Rep}_{\rm cong}(\Ga)$
coincides with the set of all finite
dimensional unitary representations of $\Ga$ 
(see \cite[Proposition 2]{Bek-RF}).
This implies (see, for instance,
\cite[Proposition 1]{Bek-RF}) that 
$$
 \bigcap_{\sigma\in{\rm Rep_{\rm cong}}(\Ga)} C^*-\Ker \sigma  \subset  C^*-\Ker \lambda_{\Ga},
$$
where $C^*-\Ker \sigma$ denotes the kernel
in $C^*(\Ga)$ of the extension of a unitary 
representation $\sigma$ of $\Ga.$ 

We consider now the restriction $\pi|_\Lambda$ of $\pi$ to $\Lambda.$
Observe that 
$${\rm Rep_{\rm cong}}(\Lambda)=\left\{\sigma|_\Lambda\ :\ 
\sigma\in{\rm Rep_{\rm cong}}(\Ga)\right\}.$$
Since $C^*-\Ker \lambda_{\Lambda}=C^*-\Ker (\lambda_{\Ga}|_{\Lambda}),$
we have
$$
 \bigcap_{\sigma\in{\rm Rep_{\rm cong}}(\Lambda)} C^*-\Ker \sigma  \subset  C^*-\Ker \lambda_{\Lambda},
$$
It follows from  Selberg's inequality $\lambda_1\geq 3/16$ (see \cite[Lemma~3]{Bek-RF})
and from the fact that $SL_2(\ZZ)$ does not have  Kazhdan's Property (T)
that ${\rm Rep_{\rm cong}}(\Lambda)$ does not separate
the points of $C^*(\Lambda),$ that is,
$$
\bigcap_{\sigma\in{\rm Rep_{\rm cong}}(\Lambda)} C^*-\Ker\sigma\neq\{0\}.
$$
Hence, we have
\begin{eqnarray*}
C^*-\Ker (\pi|_\Lambda) &=&C^*-\Ker (\pi_{\infty}|_\Lambda) \cap\bigcap_i C^*-\Ker(\sigma_{i}|_\Lambda)\\
&=&C^*-\Ker \lambda_{\Lambda} \cap\bigcap_i C^*-\Ker(\sigma_{i}|_\Lambda)\\
&\supset&
C^*-\Ker \lambda_{\Lambda}\cap \bigcap_{\sigma\in{\rm Rep_{\rm cong}}(\Lambda)} C^*-\Ker\sigma\\
&=&\bigcap_{\sigma\in{\rm Rep_{\rm cong}}(\Lambda)} C^*-\Ker\sigma
\end{eqnarray*}
and $C^*-\Ker (\pi|_\Lambda)\neq\{0\}.$
This clearly implies that $\f|_\Lambda$ is not faithful.
$\bsq$
\end{proof}
\begin{remark}
\label{Rem-Choi}
 The previous result  does not
hold for $n=2.$ Indeed, as was shown in \cite[Corollary 9]{Choi},
$C^*(SL_2(\ZZ)$ has a faithful trace.
In fact a stronger result is proved in \cite[Theorem~7]{Choi}:
$C^*(SL_2(\ZZ)$ is residually finite dimensional,
that is, the finite dimensional
representations of $SL_2(\ZZ)$ separate 
the points of $C^*(SL_2(\ZZ)).$ 

It is shown in  \cite{LubSha}
that other interesting groups  have a residually
 finite dimensional full $C^*$-algebra; this is, for instance,
the case for fundamental groups of surfaces.

\end{remark}

\section{A remark on semi-finite traces}
\label{S:Theo3}
As mentioned in the introduction,
it is conceivable that semi-finite,
infinite traces exist on $C^*(\PSL)$
for $n\geq 3.$ 
The following result implies that no such trace
factorizes through the reduced $C^*$-algebra
 $C^*_r(\PSL)$ for any integer $n\geq 2.$ 

\begin{proposition}
\label{Theo3} 
Let $G$ be a connected real semisimple Lie group
without compact factors and with trivial centre. 
Let $\Ga$ be a Zariski-dense 
subgroup of $G$. Then the tracial state $\delta_e$   is,
up to a scalar multiple, the unique 
semi-finite trace on $C_r^*(\Ga).$
In particular, $C_r^*(\Ga)$
has no normal factor representation of type $II_{\infty}.$
\end{proposition}

\begin{proof}
Let $\f: C^*_r(\Ga)^+\to [0,\infty]$ be 
a  semi-finite trace on the 
set of positive elements of $C^*_r(\Ga).$ 

We use an observation from \cite[page 583]{Rosenberg}.
It is well-known that there exist a non-zero two-sided
ideal $\frak m,$ called the ideal of definition
of $\f,$ and a linear functional
on $\frak m$ which coincides with
$\f$ on $\frak m^+$ (see \cite[Proposition~6.1.2]{Dixmier}).
Now, by \cite{BCH},
$C^*_r(\Ga)$ is simple, that is,
$C^*_r(\Ga)$ has no non-trivial two-sided
(closed or non-closed) ideals. Hence, 
 $\frak m =C^*_r(\Ga)$ and $\f$ is a finite
trace. By  \cite{BCH}, $\delta_e$ is the unique
tracial state on  $C^*_r(\Ga)$ and the claim
follows.$\bsq$
\end{proof}

Examples of Zariski dense subgroups $\Ga$ 
of a group $G$ as in the previous proposition
include all lattices in $G.$ So
Proposition~\ref{Theo3} applies, for instance, when
 $\Ga=PSL_n(\ZZ)$ for $n\geq 2$ 
or when $\Ga$ is the fundamental group
of an oriented compact surface of genus $\geq 2.$

\noindent
{\bf Address}

\noindent
 Bachir Bekka, UFR Math\'ematique, Universit\'e de  Rennes 1, 
Campus Beaulieu, F-35042  Rennes Cedex, France

\noindent
E-mail : bachir.bekka@univ-rennes1.fr


\begin{thebibliography}{BeMa00}


\bibitem[BaMS67]{BLS} H. Bass, M. Lazard and J-P. Serre.
\newblock Sous-groupes d'indice fini dans  $SL(n,\ZZ)$.
 \newblock\emph{Bull. Amer. Math. Soc.} {\bf  70}, 385-392 (1964).




\bibitem[BeCH95]{BCH} B. Bekka, M. Cowling and  P. de la Harpe.
\newblock Some groups whose reduced $C\sp*$-algebra is simple. 
 \newblock\emph{Publ. Math. IHES} {\bf 80}, 117-134 (1995).

\bibitem[Bekk99]{Bek-RF} B. Bekka.
\newblock On the full $C^*$-algebras of arithmetic groups and the congruence subgroup problem. 
\newblock\emph{Forum Math.} {\bf 11}, 705-715 (1999).


\bibitem[Boye83]{Boyer} R.P. Boyer.
\newblock Infinite traces of AF-algebras and characters of U($\infty)$. 
\newblock\emph{J. Oper. Theory} {\bf 9}, 205-236 (1983).




\bibitem[Bren60]{Brenner} J.L. Brenner.
\newblock The linear homogeneous group, III.
\newblock\emph{Ann. Math.} {\bf 71}, 210-223 (1960).






\bibitem[Burg91]{Burger}  M. Burger.
\newblock Kazhdan constants for $SL(3,\ZZ)$.
\newblock\emph{J. Reine Angew. Math.} {\bf 413}, 
36-67 (1991).


\bibitem[Choi80]{Choi} M. D. Choi.  
\newblock The full $C^*$--algebra of the free group
on two generators.
\newblock\emph{Pac. J. Math.} 
{\bf 87}, 41-48 (1980).


\bibitem[Conn76]{Connes} A. Connes. 
\newblock{Classification of injective factors.},
\newblock \emph{Ann. Math.} {\bf 104}, 73-115 (1976).

\bibitem[Conn80]{Conne--80}  A. Connes.
\newblock A factor of type $II_1$ with countable fundamental group.
\newblock \emph{J. Operator Theory}, 4: 151--153, 1980.


\bibitem[CoJo85]{CoJo} A. Connes and V. Jones. \newblock{Property (T) for von Neumann
algebras.}
\newblock \emph{Bull. London Math. Soc.} {\bf 17}, 51-62 (1985).


\bibitem[CoHa89]{Cowling-Haagerup}  M. Cowling and U. Haagerup.
\newblock Completely bounded multipliers of the Fourier algebra of a simple Lie group
of real rank one.
\newblock \emph{Invent. Math.}, 96:507--542, 1989.




\bibitem[DaKe79]{DaniKeane} S.G. Dani and M. Keane.
\newblock Ergodic invariant measures for action of $SL(2,\ZZ).$ 
\newblock \emph{Ann. Inst. Henri Poincar\'e} {\bf 15}, Nouv. S\'er., Sect. B, 79-84 (1979).


\bibitem[Dix-C*]{Dixmier} J. Dixmier.
\newblock\emph{Les $C\sp *$-alg\`ebres et leurs repr\'esentations.} 
\newblock Gauthier-Villars, 1969.

\bibitem[Dix-vN]{Dixmier-vN} J. Dixmier.
\newblock\emph{Les alg\`ebres d'op\'erateurs dans l'espace Hilbertien.} 
\newblock Gauthier-Villars, 1969.


\bibitem[Furm99]{Furman} A. Furman. 
\newblock 
Orbit  equivalence rigidity.  %
\newblock \emph{Ann. Math.}, 150:1083--1108, 1999.




\bibitem[Harp79]{Pierre} P. de la Harpe 
\newblock Simplicity of the projective unitary groups defined by simple factors.
\newblock Comment. Math. Helv.  {\bf 54},   334--345 (1979).



\bibitem[Jone00]{Jones} V. Jones. Ten problems. 
\newblock In: Mathematics: Frontiers and perspectives. 
\newblock Ed: V. Arnold et al., pages 79-91.
\newblock American Mathematical Society, 2000.




\bibitem[Kirc93]{Kirchberg}  E. Kirchberg.
 \newblock On non--split extensions, tensor products
and exactness of group $C^*$--algebras.
\newblock
\emph{Invent. Math.} {\bf 112}, 449--489 (1993).

\bibitem[Kiri65]{Kirillov}  A. A. Kirillov.
\newblock   Positive definite
functions on a group of matrices with elements from a discrete field.
\newblock\emph{Soviet. Math. Dokl.} {\bf 6}, 707-709 (1965).



\bibitem[LuSh04]{LubSha} A. Lubotzky and Y. Shalom.
 \newblock Finite representations in the unitary dual and Ramanujan groups. 
 \newblock In: ``Discrete geometric analysis",
 Proceedings of the 1st JAMS symposium, Sendai, 2002.
 \newblock  Contemporary Mathematics {\bf 347} , 173-189 (2004).



\bibitem[Marg91]{Margulis} G.A. Margulis. \newblock \emph{Discrete
subgroups of semisimple Lie groups}, 
\newblock %
Springer-Verlag, 1991.



\bibitem [Menn65]{Mennicke} J. L. Mennicke.
 \newblock Finite factor groups of the unimodular group.
\emph{Ann. Math.} {\bf 81}, 31-37 (1965).




\bibitem[Popa06-a]{Popa1} S. Popa.
\newblock   Strong rigidity of $II_1$ factors arising from malleable actions of w-rigid groups, I.
\newblock  \emph{Inventiones Math.} {\bf 165}, 369-408 (2006).


\bibitem[Popa06-b]{Popa2} S. Popa.
\newblock   Strong rigidity of $II_1$ factors arising from malleable actions of w-rigid groups, II.
\newblock  \emph{Inventiones Math.} {\bf 165}, 409-452 (2006).


\bibitem[Robe93]{Guyan} G. Robertson.
\newblock Property (T) for $II_1$ factors and unitary
representations of Kazhdan groups.
\newblock \emph{Math. Ann.}
{\bf 296}, 547-555 (1993).


\bibitem[Ovci71]{Ovc}  S.V. Ovcinnikov.
\newblock   Positive definite
functions on Chevalley groups.
\newblock\emph{Funct. Anal. Appl.} {\bf 5}, 79-80 (1971).

\bibitem[Rose89]{Rosenberg} J. Rosenberg.
\newblock Un compl\'ement \`a un th\'eor\`eme de Kirillov
sur les caract\`eres de $GL(n)$ d'un corps infini
discret.
\newblock \emph{C. R. Acad.Sc.Paris}
{\bf 309}, S\'erie I, 581-586 (1989).


\bibitem[Skud76]{Skud} H-L. Skudlarek.
\newblock Die unzerlegbaren Charaktere einiger diskreter Gruppen. 
\newblock \emph{Math. Ann.} {\bf 223}, 213-231 (1976).

\bibitem[Stei85]{Steinberg} R. Steinberg.
\newblock Some consequences of the elementary relations in $SL_n.$ 
\newblock \emph{Contempary Math.} {\bf 45}, 335-350 (1985).


\bibitem[StVo75] {StraVoi} S. Stratila and  D. Voiculescu.
\newblock \emph{Representations of $AF$-algebras and of the group $U(\infty)$.}
\newblock Lecture Notes in Mathematics {\bf 486},
 Springer, 1975. 



\bibitem[Thom64a]{ThomaTyp1} E. Thoma.
\newblock \"Uber unit\"are Darstellungen abz\"ahlbarer, diskreter Gruppen.
 \newblock \emph{Math. Ann.} {\bf 153}, 111-138 (1964).


\bibitem[Thom64b]{ThomaSym} E. Thoma.
\newblock
Die unzerlegbaren, positiv-definiten Klassenfunktionen der abz\"ahlbar unendlichen, symmetrischen Gruppe.
\newblock \emph{Math. Z.} {\bf 85}, 40-61 (1964).




\bibitem[Tits76]{Tits} J. Tits.
\newblock Syst\`emes g\'en\'erateurs de groupes de congruence. 
\newblock \emph{C. R. Acad. Sci., Paris}, S\'er. A, {\bf 283}, 693-695 (1976).


\bibitem[Vaes06]{Vaes} S. Vaes.
\newblock Rigidity results for Bernoulli actions and their von Neumann algebras (after Sorin Popa).
\newblock S\'eminaire Bourbaki, Expos\'e 961, March 2006

\bibitem[Vale97]{Alain} A. Valette.
\newblock Amenable representations and finite injective von
Neumann algebras.
\newblock \emph{Proc. Amer. Math. Soc.}
{\bf 125} , 1841-1843 (1997).


\bibitem[VerKe81]{Vershik} A.M. Vershik and S.V. Kerov.
\newblock Characters and factor representations of the infinite symmetric group. 
\newblock \emph{Sov. Math. Dokl.} {\bf 23}, 389-392 (1981).


\bibitem[Voic76]{Voiculescu} D. Voiculescu.
\newblock Repr\'esentations factorielles de type $II_1$ de $U(\infty)$.
\newblock \emph{J. Math. Pures Appl.} {\bf IX}, S\'er. 55, 1-20 (1976). 

\bibitem[Zele81]{Zel} A. V. Zelevinsky.
\newblock \emph{Representations of finite classical groups. A Hopf algebra approach.} 
\newblock Lecture Notes in Mathematics {\bf 869}, Springer, 1981.

\bibitem[Zimm84]{Zimmer} R.J. Zimmer.  
\newblock \emph{Ergodic theory
and semisimple groups}, 
\newblock %
Birkh\"auser, 1984.




\end{thebibliography}
\end{document}